%% file: PZF_ICA.tex
\documentclass[11pt, notitlepage]{article}   

\usepackage{dsfont}
\usepackage{amsmath,amsthm,amsfonts}   
\usepackage{amscd}
\usepackage{amsfonts}
\usepackage{amssymb}
\usepackage{color}      
\usepackage{epsfig}
\usepackage{graphicx}           

\theoremstyle{plain}
\newtheorem{Thm}{Theorem}[section]

\newtheorem{Prop}[Thm]{Proposition}

\theoremstyle{definition}
\newtheorem{Def}[Thm]{Definition}

\newtheorem{Eg}{Example}
\newtheorem{Que}{Question}
\theoremstyle{remark}

\def\finf{\mathop{{\rm I}\kern -.27 em {\rm F}}\nolimits}

\begin{document}

\title{Probabilistic Zero Forcing in Graphs}

\author{{\bf Cong X. Kang$^1$} and {\bf Eunjeong Yi$^2$}\\
\small Texas A\&M University at Galveston, Galveston, TX 77553, USA\\
$^1${\small\em kangc@tamug.edu}; $^2${\small\em yie@tamug.edu}}

\maketitle

\date{}

\begin{abstract}
The \emph{zero forcing number} $Z(G)$ of a graph $G$ is the minimum cardinality of a set $S$ of black vertices (whereas vertices in $V(G)\!\setminus\!S$ are colored white) such that $V(G)$ is turned black after finitely many applications of ``the (classical) color change rule": a white vertex is converted to a black vertex if it is the only white neighbor of a black vertex. Zero forcing number was introduced and used to bound the minimum rank of graphs by the ``AIM Minimum Rank -- Special Graphs Work Group". We introduce here a probabilistic color change rule (pccr) which is a natural generalization of the classical color change rule. We introduce a theory of probabilistic zero forcing arising out of the pccr; the theory yields a quantity $P_A(G)$, which can be viewed as the probability that a graph $G$ with an initial black set $A$ will be converted entirely to the color black. We also interpret the evolution of the sample spaces of this theory as a Markov process. We end with a few basic examples illustrating this theory.
\end{abstract}


\vspace{.05in}


\section{Introduction}

Let $G = (V(G),E(G))$ be a finite, simple, connected, and undirected graph of order $|V(G)|=n \ge 2$. For a vertex $v \in V(G)$, the \emph{open neighborhood of $v$} is the set $N_G(v)=\{u \mid uv \in E(G)\}$, the \emph{closed neighborhood of $v$} is the set $N_G[v]=N_G(v) \cup \{v\}$, and the \emph{degree} of a vertex $v$ is $\deg_G(v)=|N_G(v)|$; we drop $G$ if it is clear in context.\\

Let each vertex of a graph $G$ be given either the color black or the color white. Denote by $Z$ the initial set of black vertices of $G$. The ``color change rule" changes the color of a vertex $v$ from white to black if $v$ is the only white neighbor of a black vertex $u$; in this case, we say that \emph{$u$ forces $v$} and write $u \rightarrow v$. The initial black set $Z$ is said to be ``zero forcing" if, after finitely many applications of the color change rule, all vertices of $G$ are forced to black. The ``zero forcing number of $G$", $Z(G)$, is defined as the minimum of the cardinalities of all zero forcing sets of $G$. \\

The notion of a zero forcing set, as well  as the associated zero forcing number, of a graph was introduced by the ``AIM Minimum Rank -- Special Graphs Group" (see \cite{AIM}) to bound the minimum rank for numerous families of graphs. Zero forcing parameters were then studied and applied to the minimum rank problem in~\cite{preprint, UB, B2, ZFsurvey, ZFsurvey2, Huang, Pathcover}. More recently, viewing zero forcing as a dynamical process unfolding on a graph, the four authors of~\cite{iteration} and, independently, Hogben et al in~\cite{proptime} have studied the number of steps (or units of time) it takes for an initial zero forcing set to force all vertices of a graph to black. It is also noteworthy that, independent of the forgoing efforts and in connection with the control of quantum systems, this parameter has been studied by some physicists who called it ``the graph infection number"~\cite{motivation1,motivation2,motivation3}. On the other hand, probability has been playing a central role in graph theory ever since the seminal paper by Erd\"{o}s and Renyi~\cite{ErdosRenyi}. In the 2005 \emph{Nature} article entitled ``Evolutionary dynamics on graphs"~\cite{NatureArticle}, Lieberman, Hauert, and Nowak studied the question: ``what is the probability that a newly introduced mutant generates a lineage that takes over the whole population?"\\

Let us henceforth use the term ``classical zero forcing" to refer to the notion of zero forcing already defined. Taking the viewpoint that classical zero forcing is a kind of ``graph infection" or ``spread of a mutation", it's then natural -- indeed compelling -- to formulate a notion of ``probabilistic zero forcing" wherein classical zero forcing is obtained as a special case.  To this end, let $Z$ be an initial black set of a graph $G$; here $Z\neq \emptyset$, but $Z$ is not necessarily zero forcing in the classical sense. With respect to $Z$, we define a ``probabilistic color change rule" (hereafter also referred to as ``pccr") in which a black vertex forces each of its white neighbors independently to the color black with probability given as follows. Let $F(u\rightarrow v)$ denote the probability of the event $u\rightarrow v$ (which we think of as the ``forcing probability of $u$ on $v$") that the vertex $u$ forces the color of a white vertex $v$ to black. We adopt the following \\

\begin{Def}
\begin{equation*}
\displaystyle{
F(u\rightarrow v)=\left\{
\begin{array}{ll}
0 & \mbox{ if } u\notin Z \mbox{ or } v\notin N(u) \mbox{ or } N[u]\subseteq Z \\
\frac{\mid N[u]\cap Z\mid}{\deg(u)} & \mbox{ otherwise. }
\end{array} \right.
}
\end{equation*}
\end{Def}

\bigskip

We will assume that a vertex $u$ forces each of its white neighbors independently; i.e., we think of the action of a black vertex on each of its $m$ white
neighbors as a binomial (or Bernoulli) experiment with probability of success (where by a success we mean the conversion of one of the $m$ white neighbors to black) given by the formula above. Since $F(u \rightarrow v_1)=F(u \rightarrow v_2)$ when $v_1$ and $v_2$ are both white neighbors of $u$, we will simply write $F(u)$ to specify the forcing probability of $u$ when ambiguity is not a concern. Let $P( \rightarrow v)$ denote the probability (``conversion probability") that a vertex $v$ will be black after one global application (iteration) of pccr, given a set of black vertices $Z\subseteq V(G)$. We thus define $P(\rightarrow v)=1$ if $v\in Z$ and $\displaystyle P( \rightarrow v)=F(\cup_{u\in N(v)}u\rightarrow v)$ if $v\in V(G)\!\setminus\!Z$, where $F(\cup_{u\in N(v)}u\rightarrow v)=F(u_1\rightarrow v)+F(\cup_{u\in N(v)\setminus\{u_1\}}u\rightarrow v)-F(u_1\rightarrow v)\cdot F(\cup_{u\in N(v)\setminus\{u_1\}}u\rightarrow v)$; i.e., we assume that each white vertex is being acted upon independently by its black neighbors. \\

Observe that classical zero forcing where a black vertex $u$ forces its sole white neighbor $v$ is
recovered as a pair, a black vertex $u$ and a white vertex $v$, with forcing probability $F(u\rightarrow v)=1$ and conversion probability $P(\rightarrow v)=1$. From the
perspective of modeling a ``real world" situation, our formulation is sensible: it only takes one affected entity to spread a condition to an entire
population. We want to mention here that Curtis Miller has independently formulated a number of probabilistic variants of zero forcing~\cite{curtis}. But our present definition, as well as our approach and emphasis, is distinct from the ones he considered (personal communication). Notice that, once a classical zero forcing (black) set $Z_1$ containing an initial black set $Z_0$ of $G$ is attained through probabilistic zero forcing, classical zero forcing will be in force as a part of probabilistic zero forcing; before then, some combination of classical and probabilistic zero forcing occurs. In the rest of this paper, we will at times abbreviate probabilistic zero forcing as ``PZF".

\section {Interpretations of PZF}

Let's call a graph $G$, with vertices labeled $1$ through $n$, a colored graph if each of its vertices receives the color black or the color white, and let's regard as distinct two colored graphs $G_1$ and $G_2$ with the same underlying graph $G$ when there is a vertex $v\in\{1,\ldots, n\}$ which has one color in $G_1$ but has the other color in $G_2$. For a colored graph $G$, denote by $D(G)$ the set of all colored graphs which can be obtained with positive probability from $G$ in one global application of pccr. Let a graph $G^*$, with initial black set $Z$, be given. Then PZF generates a countable collection of probability sample spaces $\mathfrak{C}(G_0)=\{S^k: k\in\mathbb{N}\}$, one sample space $S^k$ at each step $k$ which can be inductively defined as follows: $S^0=\{G_0\}$, where $G_0$ is $G^*$ with $Z\subseteq V(G^*)$ colored black and $V(G^*)\setminus Z$ colored white; once $S^{k-1}$ is defined, $S^k= \bigcup_{G'\in S^{k-1}}D(G')$. Let $P^{(k)}$ be the probability function defined on $S^k$ which is derived from the pccr. It follows that $\sum_{G\in S^k}P^{(k)}(G)=1$ for each $k\in \mathbb{N}$. The $\mathfrak{C}(G_0)$ just defined represents a discrete dynamical system on the underlying graph $G$, and its evolution is the central problem of the PZF theory.\\

Let $W\subseteq V(G^*)$ and $S^k_W=\{G\in S^k: \mbox{each } v\in W \mbox{ is a black vertex of }G\}$; if $W=\{v\}$, we'll write $S^k_v$. Define $P^{(k)}(W)=\sum_{G\in S^k_W}P^{(k)}(G)$. If $W_1, W_2 \subseteq V(G^*)$, then it follows that $P^{(k)}((\rightarrow W_1) \cup (\rightarrow W_2))=P^{(k)}(\rightarrow W_1)+P^{(k)}(\rightarrow W_2)-P^{(k)}(\rightarrow (W_1\cup W_2))$, since $\rightarrow (W_1 \cup W_2)=(\rightarrow W_1) \cap (\rightarrow W_2)$. It's clear that, relative to an initial black set $Z$, $P^{(1)}(\rightarrow v)=P(\rightarrow v)$, where $P(\rightarrow v)$ is previously defined. Also notice that $P^{(k+1)}(\rightarrow v)\geq P^{(k)}(\rightarrow v)$ by \emph{color change asymmetry}; namely, if a vertex $v$ is black in a colored graph $G$, then $v$ is black for each member of $D(G)$. Suppose for some $k_0$, there is no $G\in S^{k_0}$ such that $G$ contains a black vertex $v$ which forces with probability one (i.e., $v$ has a lone white neighbor), then $S^k\subseteq S^{k+1}$ (as sets) for $k\geq k_0$. Since there are at most $2^{|V(G^*)|}$ colored graphs on the underlying graph $G^*$, the sequence $S^k$ is eventually constant; i.e., there is some $M\in\mathds{N}$ such that $S^{k+1}=S^k$ for any $k\geq M$. Note that this does not mean that $P^{(k+1)}=P^{(k)}$ for $k\geq M$. \\

The evolution of our discrete dynamical system is a Markov process. Here is the gross view: Let $G_1, G_2, \ldots, G_{2^n}$ be an ordered listing of the $2^n$ possible colored graphs (possible states of the system) underlied by a graph $G$ of order $n$: each $G_i$ can be identified with the positive integer $i$. For each $k\in \mathbb{N}$, we associate a \emph{probability vector} $\mathbf{x}^{(k)}$ of length $2^n$ whose $i$-th entry denotes the probability that our system is in the state $i$ at time $k$. The stationary transition matrix $A$ for the system has size $2^n\!\times\!2^n$, and each of its entries $A_{i,j}$ provides the probability that a system in state $j$ at time $k$ will transition into state $i$ at time $k+1$; $A_{i,j}$ is thus determined by the pccr applied to the state $j$. \\

We observe that, starting with a colored graph $G_0$ of order $n$ in $S^0$, by the fact that the black set is expanding (or non-contracting) from any $G\in S^{k-1}$ to each $G'\in D(G)$, the number of colored graphs (and hence the size of the transition matrix) requiring consideration is typically and tremendously reduced from $2^n$. Recall that, in the theory of Markov processes, an \emph{absorbing state} is one that can not be exited once entered; i.e., it is a state $i$ such that $A_{i,i}=1$; a state $i$ is \emph{transient} (also called \emph{non-recurrent}) if there is a state $j$ reachable from $i$, but state $i$ is not reachable from $j$. If there is only one absorbing state $j_0$, then $\displaystyle\lim_{k\rightarrow \infty} A^{k}_{i,j_0}=1$ for each state $i$ (see p.247 and p.241 of~\cite{sto1}, also see pp.58-66 of~\cite{sto2} for proofs). Since our dynamical system has transient states except for one absorbing state (the colored graph $G_b$ where all vertices of $G_b$ are black), we have the following

\begin{Prop}\label{limiting state}
For any colored graph $G_0$, $\displaystyle\lim_{k\rightarrow \infty}P^{(k)}(G_b)=1$, where $G_b\in S^k\in \mathfrak{C}(G_0)$.
\end{Prop}

That is, a colored graph will turn entirely black with probability one as the number of steps $k$ goes to infinity. The following is an instructive example.\\

\Eg Consider the complete bipartite graph $K_{1,m}$, where $m\geq 1$, which we'll also call the star with $m$ pendants; let the center $v_0$ be the only black vertex of $G_0$. We have $S^0=\{G_0\}$ and $S^k=\{G_0,\ldots, G_{2^m-1}\}$ for $k\geq 1$; i.e., all possible colored graphs (or states) occur in every sample space $S^k$ for $k\geq 1$. Let $G_b$ denote the state where all vertices of the star are black. Then Proposition~\ref{limiting state} says that $\lim_{k\rightarrow\infty} P^{(k)}(G_b)=1$, which can also be seen via a direct, inductive argument on the number of pendants $m$ via the following observation: the center $v_0$ is the only vertex capable of forcing; each $G'\in D(G)\subseteq S^{k+1}$ is obtained from $G\in S^k$ through a binomial experiment wherein each white pendant is turned black by the forcing strength of $v_0$; the average forcing strength of $v_0$ is monotonically increasing toward $1$ as $k$ goes to infinity. \hfill $\Box$\\

In order that probabilistic zero forcing be a natural extension of classical zero forcing, there ought to be a probability $P_A(G)$ assigned to each pair $(G, A)$ which indicates the likelihood that the graph $G$ will be forced entirely black given the initial black set $A\subseteq V(G)$ under finitely many applications of the pccr; moreover, this probability should satisfy the following conditions: 1) $P_\emptyset(G)=0$, 2) $P_A(G)=1$ if and only $A$ is a zero forcing set of $G$, 3) $P_A(G)\leq P_B(G)$ for $A\subseteq B\subseteq V(G)$. Given a pair $(G,A)$ and the associated $\mathfrak{C}(G_0)$, define $T^k$ to be the subset of $S^k$ such that $G'\in T^k$ if and only if the black set of $G'$ contains a classical zero forcing set of $G$; note that $T^k\neq \emptyset$ for some $k\in\mathds{N}$ by the definition of pccr. It's clear that the following definition of $P_A(G)$ meets the three requirements.\\

\Def Let $P_A(G)=P^{k_0}(T^{k_0})$, where $k_0$ is the least $k$ such that $T^k\neq \emptyset$. We think of $P_A(G)$ as the probability that a graph $G$ with an initial black set $A\subseteq V(G)$ ``will turn entirely black in finite time".\\

\Eg Let $G=C_n$, $n\geq 3$, with initial black set $A=\{v_1\}\subseteq V(G)$ be given; also let $E(G)=\{v_1v_2, v_2v_3\ldots, v_{n-1}v_{n},v_{n}v_{1}\}$. Recall two adjacent vertices form a classical zero forcing set for $G$. Then, in going from $S^0$ to $S^1$, the only forcing action is given by $F(v_1\rightarrow v_2)=\frac{1}{2}=F(v_1\rightarrow v_n)$. Only one (call it $G_w$) of the four colored graphs of $S^1$ does not contain a classical zero forcing set, and thus $P_A(G)=1-P(G_w)=1-\frac{1}{2}\cdot\frac{1}{2}=\frac{3}{4}$. Note the equivalence, trivially here, of $P_A(G)$ to $P^{(1)}((\rightarrow v_2) \cup (\rightarrow v_n))=P^{(1)}(\rightarrow v_2)+P^{(1)}(\rightarrow v_n)-P^{(1)}(\rightarrow \{v_2,v_n\})$. \hfill $\Box$\\

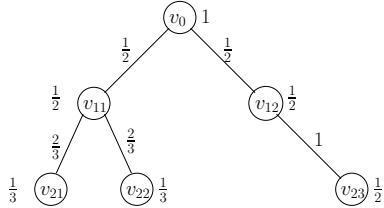
\begin{figure}[htbp]
\begin{center}
\scalebox{0.45}{\input{probzT.pstex_t}} \caption{$G$ with initial black set $\{v_0\}$}\label{tree}
\end{center}
\end{figure}

\Eg Let $G$ be the tree given in Fig~\ref{tree}, with initial black set $A=\{v_0\}$. The forcing action propagates from $v_0$ to the three pendants (vertices of degree $1$). In Fig~\ref{tree}, the number next to each edge $uv$ is $F(u\rightarrow v)$; notice that forcing can only occur in one direction for each edge in this pair $(G,A)$. For example, the number $\frac{2}{3}$ is next to the edge $v_{11}v_{21}$ because $F(v_{11}\rightarrow v_{21})=\frac{2}{3}$ given that the vertex $v_{11}$ is already (forced to) black. Likewise, the number next to each vertex indicates the probability that each vertex will be forced black given the initial black vertex $v_0$. \\

\begin{figure}[htbp]
\begin{center}
\scalebox{0.45}{\input{samplespaces2.pstex_t}} \caption{$S^0$, $S^1$, and $T^2 \subset S^2$ associated with the pair $(G, A)$}\label{tree-samplespaces}
\end{center}
\end{figure}
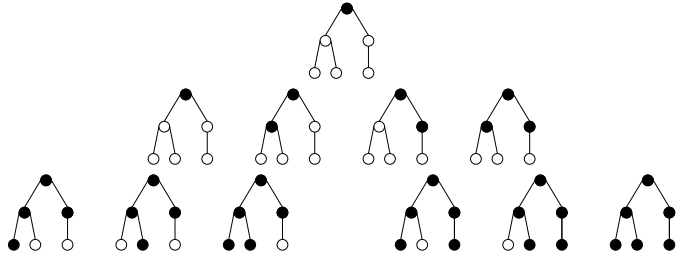

We determine $P_A(G)$ by starting with the observation that a colored graph containing a classical zero forcing set first occurs in $S^2$, since a black set $B$ containing $v_0$ is classically zero forcing if and only if $B\cap\{v_{21},v_{22}\}\neq \emptyset$. The top two rows of Figure~\ref{tree-samplespaces} show $S^0$ and $S^1$, whereas the third row only shows $T^2\subset S^2$, namely the colored graphs of $S^2$ which contain classical zero forcing sets. From left to right, the first three graphs of the third row form $D(G_2)$ and last three graphs of the third row form $D(G_4)$, where $G_2$ and $G_4$ are the second and the fourth graphs, respectively, of the second row. Now, it's clear that $P^{(1)}(G_2)=\frac{1}{4}=P^{(1)}(G_4)$.
Thus, $P_A(G)=P^{(2)}(\bigcup_{G' \in T^2}G')=2[\frac{1}{4}(\frac{2}{3}\frac{1}{3}+\frac{1}{3}\frac{2}{3}+\frac{2}{3}\frac{2}{3})]=\frac{4}{9}$ \\

We can also deduce $P_A(G)$ in this case without explicitly enumerating elements of the sample spaces $\mathfrak{C}(G_0)=\{S^k: k\in\mathbb{N}\}$ in the following way. First, still notice that a black set $B$ containing $v_0$ is classically zero forcing if and only if $B\cap\{v_{21},v_{22}\}\neq \emptyset$; also notice that $P^{(1)}(\rightarrow v_{21})=0=P^{(1)}(\rightarrow v_{22})$, but $P^{(2)}(\rightarrow v_{21})\neq 0$ and $P^{(2)}(\rightarrow v_{22})\neq 0$. Thus, the probability that $G$ will be turned entirely black by $\{v_0\}$ under PZF is given by $P^{(2)}((\rightarrow v_{21}) \cup (\rightarrow v_{22}))$, which must be the same as the probability $P(E_2)$ of the event $E_2$ that $v_{11}\rightarrow v_{21}$ or $v_{11}\rightarrow v_{22}$. Let $E_1$ denote the event that $v_0\rightarrow v_{11}$. Observe $E_1\subseteq E_2$; i.e., if $E_2$ happens, then $E_1$ already happened. Thus, we have $P(E_2)=P(E_2\cap E_1)=P(E_2|E_1)P(E_1)= (F(v_{11}\!\rightarrow\! v_{21})+F(v_{11}\!\rightarrow\! v_{22})-F(v_{11}\!\rightarrow\! v_{21})\!\cdot\! F(v_{11}\!\rightarrow\! v_{21}))\!\cdot\! F(v_0\!\rightarrow\! v_{11})=(\frac{2}{3}+\frac{2}{3}-\frac{2}{3}\cdot\frac{2}{3})\cdot\frac{1}{2}=\frac{4}{9}$. \hfill $\Box$\\

The following is a basic question in our PZF theory at this juncture -- one which we'll solve for the graph depicted in Figure~\ref{tree} in the case the initial black set is a singleton.\\

\Que\label{the question} Let $P_{(j)}(G)=\max(\{P_A(G): |A|=j\})$, where $A$ is an initial black set of $G$ and $j$ is an integer between $1$ and $Z(G)$. What is $P_{(j)}(G)$ for each $j\in\{1,\ldots, Z(G)\}$ and what are the initial black sets $A$ attaining $P_{(j)}(G)$ for each $j$? \\

\Eg Let's answer Question~\ref{the question} for $j=1$ when $G$ is the star $K_{1,m}$ with $m$ pendants, and let $m\geq 3$ for simplicity. Label the center of the star $v_0$ and the pendants $v_1$ through $v_m$. Given the symmetry group of the star, it's enough to consider two cases: 1) $A=\{v_0\}$, 2) $A=\{v_1\}$. In case 1), $F(v_0)=\frac{1}{m}$ and $T^1\neq \emptyset$; thus, $P_{\{v_0\}}(G)=P^{(1)}(T^1)={m \choose m-1}(\frac{1}{m})^{m-1}(\frac{m-1}{m})=\frac{m-1}{m^{m-1}}$. In case 2), we have $F(v_1\rightarrow v_0)=1$ and $T^1=\emptyset$ since $m\geq 3$ by hypothesis. In passing from $S^1$ to $S^2$, we have $F(v_0)=\frac{2}{m}$ and $T^2\neq \emptyset$; thus, $P_{\{v_1\}}(G)=P^{(2)}(T^2)={m-1 \choose m-2}(\frac{2}{m})^{m-2}(\frac{m-2}{m})=\frac{(m-1)(m-2)2^{m-2}}{m^{m-1}}$. Obviously, $P_{(1)}(G)=P_{\{v_1\}}(G)>P_{\{v_0\}}(G)$ when $m\geq 3$. We note that it's a straightforward matter to determine $P_{(j)}(G)$ for each $1\leq j<Z(G)=m-1$ since, for each $A$ of cardinality $j$, we only need to consider the case $v_0\in A$ and the case $v_0 \notin A$.   \hfill $\Box$ \\

\Eg Returning to the graph depicted in Figure~\ref{tree}, let's consider all cases $A_i$ where $|A_i|=1$. There are four distinct cases to consider:
1) $A_1=\{v_0\}$, 2) $A_2=\{v_{23}\}$, 3) $A_3=\{v_{11}\}$, 4) $A_4=\{v_{21}\}$. Having found in case 1) that $P_{A_1}(G)=\frac{4}{9}$ in the previous example, let's consider case 2). By the pccr, we have $P(\rightarrow v_{23})=P^{(1)}(\rightarrow v_{12})=P^{(2)}(\rightarrow v_{0})=P^{(3)}(\rightarrow v_{11})=1$, and we have in $S^4$ that $F(v_{11}\rightarrow v_{21})=\frac{2}{3}=F(v_{11}\rightarrow v_{22})$; we get $P_{A_2}(G)=P^{(4)}(T^4)=P^{(4)}((\rightarrow v_{21}) \cup (\rightarrow v_{22}))=\frac{2}{3}+\frac{2}{3}-\frac{2}{3}\cdot\frac{2}{3}=\frac{8}{9}$. In case 3), $T^1\neq \emptyset$ and $P_{A_3}(G)=P^1(T^1)=3\cdot(\frac{1}{3}\cdot\frac{1}{3}\cdot\frac{2}{3})+\frac{1}{3}\cdot\frac{1}{3}\cdot\frac{1}{3}=\frac{7}{27}$. In case 4), since $F(v_{21}\rightarrow v_{11})=1$, $S^1$ contains only the colored graph where the only black vertices are $v_{21}$ and $v_{11}$, and $T^1=\emptyset$. Since $T^2\neq \emptyset$, $P_{A_4}(G)=P^{(2)}(T^2)=1-\frac{1}{3}\cdot\frac{1}{3}=\frac{8}{9}$. Therefore, we conclude that $P_{(1)}(G)=\max(\{P_A(G): |A|=1\})=\frac{8}{9}.$ \hfill $\Box$\\

\textit{Acknowledgement.} The authors thank Kiran Chilakamarri for a number of very helpful discussions. In particular, they thank him for explaining to them the idea of iterated sample spaces and for locating an informal reference on the theory of absorbing Markov chains. The authors also thank Curtis Miller for the aforementioned preprint, which gave additional impetus for them to finish this paper. Lastly, the authors thank the anonymous referee for some helpful comments and corrections.

\end{document}

%% file: probzT.pstex_t
\begin{picture}(0,0)%
\includegraphics{probzT.pstex}%
\end{picture}%
\setlength{\unitlength}{3947sp}%
\begingroup\makeatletter\ifx\SetFigFont\undefined%
\gdef\SetFigFont#1#2#3#4#5{%
  \reset@font\fontsize{#1}{#2pt}%
  \fontfamily{#3}\fontseries{#4}\fontshape{#5}%
  \selectfont}%
\fi\endgroup%
\begin{picture}(5130,2866)(4186,-2994)
\put(8851,-2836){\makebox(0,0)[lb]{\smash{{\SetFigFont{17}{20.4}{\rmdefault}{\mddefault}{\updefault}{\color[rgb]{0,0,0}$v_{23}$}%
}}}}
\put(6901,-436){\makebox(0,0)[lb]{\smash{{\SetFigFont{17}{20.4}{\rmdefault}{\mddefault}{\updefault}{\color[rgb]{0,0,0}1}%
}}}}
\put(7201,-886){\makebox(0,0)[lb]{\smash{{\SetFigFont{17}{20.4}{\rmdefault}{\mddefault}{\updefault}{\color[rgb]{0,0,0}$\frac{1}{2}$}%
}}}}
\put(8101,-1561){\makebox(0,0)[lb]{\smash{{\SetFigFont{17}{20.4}{\rmdefault}{\mddefault}{\updefault}{\color[rgb]{0,0,0}$\frac{1}{2}$}%
}}}}
\put(8476,-2161){\makebox(0,0)[lb]{\smash{{\SetFigFont{17}{20.4}{\rmdefault}{\mddefault}{\updefault}{\color[rgb]{0,0,0}$1$}%
}}}}
\put(9301,-2836){\makebox(0,0)[lb]{\smash{{\SetFigFont{17}{20.4}{\rmdefault}{\mddefault}{\updefault}{\color[rgb]{0,0,0}$\frac{1}{2}$}%
}}}}
\put(7651,-1636){\makebox(0,0)[lb]{\smash{{\SetFigFont{17}{20.4}{\rmdefault}{\mddefault}{\updefault}{\color[rgb]{0,0,0}$v_{12}$}%
}}}}
\put(6451,-436){\makebox(0,0)[lb]{\smash{{\SetFigFont{17}{20.4}{\rmdefault}{\mddefault}{\updefault}{\color[rgb]{0,0,0}$v_0$}%
}}}}
\put(6301,-2836){\makebox(0,0)[lb]{\smash{{\SetFigFont{17}{20.4}{\rmdefault}{\mddefault}{\updefault}{\color[rgb]{0,0,0}$\frac{1}{3}$}%
}}}}
\put(5251,-1636){\makebox(0,0)[lb]{\smash{{\SetFigFont{17}{20.4}{\rmdefault}{\mddefault}{\updefault}{\color[rgb]{0,0,0}$v_{11}$}%
}}}}
\put(5851,-2836){\makebox(0,0)[lb]{\smash{{\SetFigFont{17}{20.4}{\rmdefault}{\mddefault}{\updefault}{\color[rgb]{0,0,0}$v_{22}$}%
}}}}
\put(4651,-2836){\makebox(0,0)[lb]{\smash{{\SetFigFont{17}{20.4}{\rmdefault}{\mddefault}{\updefault}{\color[rgb]{0,0,0}$v_{21}$}%
}}}}
\put(5851,-2161){\makebox(0,0)[lb]{\smash{{\SetFigFont{17}{20.4}{\rmdefault}{\mddefault}{\updefault}{\color[rgb]{0,0,0}$\frac{2}{3}$}%
}}}}
\put(5776,-886){\makebox(0,0)[lb]{\smash{{\SetFigFont{17}{20.4}{\rmdefault}{\mddefault}{\updefault}{\color[rgb]{0,0,0}$\frac{1}{2}$}%
}}}}
\put(4201,-2836){\makebox(0,0)[lb]{\smash{{\SetFigFont{17}{20.4}{\rmdefault}{\mddefault}{\updefault}{\color[rgb]{0,0,0}$\frac{1}{3}$}%
}}}}
\put(4801,-1561){\makebox(0,0)[lb]{\smash{{\SetFigFont{17}{20.4}{\rmdefault}{\mddefault}{\updefault}{\color[rgb]{0,0,0}$\frac{1}{2}$}%
}}}}
\put(4801,-2236){\makebox(0,0)[lb]{\smash{{\SetFigFont{17}{20.4}{\rmdefault}{\mddefault}{\updefault}{\color[rgb]{0,0,0}$\frac{2}{3}$}%
}}}}
\end{picture}%

%% file: samplespaces2.pstex_t
\begin{picture}(0,0)%
\includegraphics{samplespaces2.pstex}%
\end{picture}%
\setlength{\unitlength}{3947sp}%
\begingroup\makeatletter\ifx\SetFigFont\undefined%
\gdef\SetFigFont#1#2#3#4#5{%
  \reset@font\fontsize{#1}{#2pt}%
  \fontfamily{#3}\fontseries{#4}\fontshape{#5}%
  \selectfont}%
\fi\endgroup%
\begin{picture}(9316,3466)(1868,-3744)
\end{picture}%